\documentclass[preprint,12pt]{elsarticle}
 \journal{}
\usepackage {multicol}
\usepackage{amssymb}
\usepackage{amsmath}
\usepackage{mathrsfs}
\usepackage{hyperref}
\usepackage{graphicx}
\usepackage{color}
\usepackage{caption}
\usepackage{esint}

\begin{document}
\begin{frontmatter}

\title{Extension of Limit Theory with Deleting Items Partial Sum of Random Variable Sequence}

\author{Jingwei Liu \corref{cor1}}
\ead{liujingwei03@tsignhua.org.cn}
\cortext[cor1]{Corresponding author.}
\address{School of Mathematics and System Sciences,Beihang University,Beijing,100083,P.R China}

\begin{abstract}

The deleting items theorems of weak law of large numbers (WLLN),strong law of large numbers (SLLN) and central limit theorem (CLT) are derived by substituting partial sum of random variable sequence with deleting items partial sum. We address the background of deleting items limit theory of random variable sequence, discuss the classical limit theory of Chebyshev WLLN, Bernoulli WLLN and Khinchine WLLN with standard mathematical analytical technique, then develop  the deleting items theorems of WLLN, SLLN and CLT based on convergence theorems and Slutsky's theorem. Our theorems extend the classical limit theory of random variable sequence and provide the construction of some asymptotic bias estimators of sample expectation and variance.

\end{abstract}


\end{frontmatter}

\section{Introduction}
\label{}

Limit theory of random variable sequence, including weak law of large numbers (WLLN), strong law of large numbers (SLLN) and central limit theorem (CLT), takes important roles in probability theory, since it provides rigorous foundation of probability theory in modern mathematics analysis [1-13].

Our motivation of deleting items limit theory in this paper is based on the practical problems and theoretical consideration.

Firstly, we start from some more general problems in random experiments or practical problems, for example,

 Example 1, When tossing a coin, if we do the experiments on a table, occasionally, the coin may jump off the table. Will the random experiment be meaningful if it goes on?

 Example 2, While simulating the Buffon needling problem, in some extreme cases, the needle does not lie across the two strip lines£¬it may be far from the area. How to deal with these cases?

 Example 3, In recent high accurate industrial techniques, some samples for parameter estimation may be measured under complicated environment or conditions, it is not easy to guarantee each sample strictly follow the theory condition. Can we believe the experimental results?

 Example 4, In modern statistical learning and data analysis science, how many samples are enough to estimate the statistical parameters or describe a statistical model? Can the subset (not subsequence) of samples perform the task?

Suppose $\xi_1,\xi_2,\cdots,\xi_n,\cdots$ is a real-valued random variable sequence in $\mathbf{R}$. We only address $\mathbf{R}$ case in this paper. Fix a natural number $k$ $(0<k<n)$ or a natural number $k(n)$ $(0<k(n)<n)$, denote $J=\{1,2,\cdots,n\}$, $J_k=\{i_1,i_2,\cdots,i_k\}$,  $J_{k(n)}=\{i_1,i_2,\cdots,i_{k(n)}\}$. where $\{i_1,i_2,\cdots,i_k\}$ or $\{i_1,i_2,\cdots,i_{k(n)}\}$ is any $k$ (or $k(n)$) different elements of $J$.

Denote\\
\begin{equation}
S_n=S_J=\sum\limits_{i=1}^n \xi_i, \ \
S_{J\backslash J_k}=\sum\limits_{i\in J\backslash J_k} \xi_i , \ \
S_{J\backslash J_{{k(n)}}}=\sum\limits_{i\in J\backslash J_{k(n)}} \xi_i .
\end{equation}

If $k$ and $k(n)$ denote the number of extreme cases in above examples, to be different from the missing completely at random (MCAR) in statistics [14-15] , we call $\displaystyle S_{J\backslash J_k}$ and $\displaystyle S_{J\backslash J_{{k(n)}}}$ deleting items  partial sum [16] of random variable sequence.

Secondly, the classical limit theory is to discuss the WLLN, SLLN and CLT of $\displaystyle \frac{S_n-a_n}{B_n}$, where $a_n$ and $B_n$ ($B_n>0$) are real-valued sequences.  If we know the probability or statistical structure of  $\displaystyle \frac{S_n-a_n}{B_n}$, then a theoretical consideration is how about $\displaystyle \frac{S_{J\backslash J_k}-a_n}{B_n}$ and $\displaystyle \frac{S_{J\backslash J_{k(n)}}-a_n}{B_n}$ ?

As $\displaystyle S_{J\backslash J_k}$ and $\displaystyle S_{J\backslash J_{{k(n)}}}$ are deleting item sum of $S_n$ or $S_J$. $\displaystyle \frac{1}{n}S_{J\backslash J_k}$ and $\displaystyle \frac{1}{n} S_{J\backslash J_{k(n)}}$ are deleting item sum of $\displaystyle \frac{1}{n}S_n$, and they are obviously not the subsequence of $\displaystyle \frac{1}{n}S_n$. Therefore, the limit theory of $\displaystyle \frac{1}{n}S_{J\backslash J_k}$ and $\displaystyle \frac{1}{n} S_{J\backslash J_{k(n)}}$ are an interesting problem, which will extend massive statistical theorems based on the limit theory of $\displaystyle \frac{1}{n}S_n$.

The rest of the paper are organized as follows. Section 2 introduces the notation and convergence theorem of random sequences. Section 3 investigates the WLLN of tradition WLLNs with classical mathematical proof. Section 4 discusses the uniform deleting items WLLN. Section 5 discusses the deleting items SLLN.Section 6 discusses the deleting items and CLT. Section 7 discusses the application of  asymptotic statistics of deleting items WLLN. Conclusion and discuss are in Section 8.

\section{Notation and Convergence}
\label{}
The definitions and relationship of convergence almost surely, convergence in probability, convergence in $r-$means and convergence in distribution are from [1-13].

Definition 1. Let $\xi_1,\xi_2,\cdots $ be random variables. $\xi_n$ converges almost surely (a.s.) to the random variable $\xi$ as $n\rightarrow +\infty$ iff
\begin{equation}
 P(\{\omega: \lim \limits_{n\rightarrow +\infty} \xi_n (\omega)=\xi (\omega) \})=1.
\end{equation}
Denote as $\xi_n \stackrel{a.s}{\longrightarrow} \xi, n\rightarrow +\infty$

Definition 2. Let $\xi_1,\xi_2,\cdots $ be random variables. $\xi_n$ converges in probability to the random variable $\xi$ iff, for any $\varepsilon>0$,
\begin{equation}
 \lim \limits_{n\rightarrow +\infty} P(|\xi_n -\xi |>\varepsilon)=0 .
\end{equation}
Denote as $\xi_n \stackrel{P}{\longrightarrow} \xi, n\rightarrow +\infty$.

Definition 3. Let $\xi_1,\xi_2,\cdots $ be random variables. $\xi_n$ converges in $r-$mean to the random variable $\xi$ as $n\rightarrow +\infty$ iff
\begin{equation}
 \lim \limits_{n\rightarrow +\infty} E(|\xi_n -\xi |^{r})=0.
\end{equation}
where $r>0$. Denote as $\xi_n \stackrel{r}{\longrightarrow} \xi, n\rightarrow +\infty$. When $r=2$, it is called convergence in square mean.

Definition 4. Let $\xi_1,\xi_2,\cdots $ be random variables. $\xi_n$ converges in distribution to the random variable $\xi$ as $n\rightarrow +\infty$ iff
\begin{equation}
 \lim \limits_{n\rightarrow +\infty} F_{\xi_n}(x) =F_{\xi}(x). \quad  \forall x\in C(F_{\xi}).
\end{equation}
where $C(F_{\xi})=\{x: F_{\xi}(x) \mbox{is continous at} \ x \}$. Denote as $\xi_n \stackrel{d}{\longrightarrow} \xi, n\rightarrow +\infty$.

Lemma 1.  Let \{$\xi_n$\} be a sequence of random variables. The following implication holds.

\begin{equation}
\begin{array}{ll}
 \xi_n \stackrel{a.s}{\longrightarrow} \xi \Rightarrow \xi_n \stackrel{P}{\longrightarrow} \xi  \Rightarrow \xi_n \stackrel{d}{\longrightarrow} \xi \\
 \xi_n \stackrel{r}{\longrightarrow} \xi \Rightarrow \xi_n \stackrel{P}{\longrightarrow} \xi.
\end{array}
\end{equation}

Lemma 2.  Let $\{\xi_n\}$£¬$\{\eta_n\}$ be  sequences of random variables. And,  $\xi_n\stackrel{a.s}{\longrightarrow}\xi$ , $\eta_n\stackrel{a.s}{\longrightarrow}\eta$. Then£¬

\begin{equation}
\begin{array}{ll}
 \xi_n \pm \eta_n \stackrel{a.s}{\longrightarrow} \xi \pm \eta \\
 \xi_n \eta_n \stackrel{a.s}{\longrightarrow} \xi \eta \\
 \xi_n /\eta_n \stackrel{a.s}{\longrightarrow} \xi/ \eta \\
\end{array}
\end{equation}

Lemma 3.  Let $\{\xi_n\}$£¬$\{\eta_n\}$ be  sequences of random variables. And,  $\xi_n\stackrel{P}{\longrightarrow}\xi$ , $\eta_n\stackrel{P}{\longrightarrow}\eta$. Then£¬

\begin{equation}
\begin{array}{ll}
 \xi_n \pm \eta_n \stackrel{P}{\longrightarrow} \xi \pm \eta,\\
 \xi_n \eta_n \stackrel{P}{\longrightarrow} \xi \eta, \\
 \xi_n /\eta_n \stackrel{P}{\longrightarrow} \xi/ \eta.
\end{array}
\end{equation}

Lemma 4. Let \{$\xi_n$\},\{$\eta_n$\} be sequences of random variables. If $\xi_n$ converges in distribution to a random variable $\xi$, and $\eta_n$ converges in probability to a constant $c$, then

\begin{equation}
\begin{array}{ll}
\displaystyle  \xi_n+\eta_n \stackrel{d}{\longrightarrow} \xi+c,  \\
\displaystyle  \xi_n \eta_n \stackrel{d}{\longrightarrow} c\xi,  \\
\displaystyle  \xi_n/\eta_n \stackrel{d}{\longrightarrow} \xi/c, \ \  c\neq 0.  \\
\end{array}
\end{equation}
where $\stackrel{d}{\longrightarrow}$ denotes convergence in distribution.

To easily compare our extension limit theory to classical limit theory, the theorems of WLLN, SLLN and CLT list as follows are mainly referred [1-13] .

\textbf{Theorem 1} (Bernoulli WLLN) Let $\xi_{1}, \xi_{2}, \cdots, \xi_{n}, \cdots $ be mutual independent random sequence, where $E\xi_{i}=p$, $D\xi_{i}=p(1-p)$, $0<p<1$ . Then, for $\forall \varepsilon>0$,
\begin{equation}
\begin{array}{ll}
\displaystyle  \lim\limits_{n \longrightarrow +\infty } P\left( \left|\frac{1}{n}\sum\limits_{i=1}^n \xi_i -p \right|< \varepsilon \right) = 1.
\end{array}
\end{equation}

\textbf{Theorem 2} (Chebyshev WLLN) Let $\xi_{1}, \xi_{2}, \cdots, \xi_{n}, \cdots $ be piecewise uncorrelated random sequence with bounded variance, where $ D\xi_{i}\leq C$, $i=1,2,\cdots $, $C$ is a constant. Then, for $\forall \varepsilon>0$,
\begin{equation}
\begin{array}{ll}
\displaystyle  \lim\limits_{n \rightarrow +\infty } P\left( \left|\frac{1}{n} \sum\limits_{i=1}^n \xi_i - \frac{1}{n} \sum\limits_{i=1}^n E\xi_i \right|< \varepsilon \right) = 1.
\end{array}
\end{equation}

\textbf{Theorem 3} (Khinchine WLLN) Let $\xi_{1}, \xi_{2}, \cdots, \xi_{n}, \cdots $ be independent identical distribution random variable sequence, with $E\xi_{i}=a$ $(i=1,2,\cdots)$. Then,
\begin{equation}
\begin{array}{ll}
\displaystyle \frac{1}{n} \sum\limits_{i=1}^n \xi_i \stackrel{P}{\longrightarrow} a , \ \ n \longrightarrow +\infty.
\end{array}
\end{equation}

\textbf{Theorem 4} Let $\xi_{1}, \xi_{2}, \cdots, \xi_{n}, \cdots $ be a random variable sequence. Then,
\begin{equation}
\begin{array}{ll}
\displaystyle    \frac{1}{n} \sum\limits_{i=1}^n \xi_i-\frac{1}{n} \sum\limits_{i=1}^n E\xi_i  \stackrel{P}{\longrightarrow} 0 , \ \ n \longrightarrow +\infty.\\
\Longleftrightarrow  \displaystyle E[\frac{|\frac{1}{n} \sum\limits_{i=1}^n \xi_i-\frac{1}{n} \sum\limits_{i=1}^n E\xi_i|^2}{1+|\frac{1}{n} \sum\limits_{i=1}^n \xi_i-\frac{1}{n} \sum\limits_{i=1}^n E\xi_i|^2}  ] \stackrel{P}{\longrightarrow} 0 , \ \ n \longrightarrow +\infty.
\end{array}
\end{equation}

\textbf{Theorem 5} (SLLN). Let $\xi_1,\xi_2,\cdots$ be pairwise independent identically distributed ( i.i.d.) random variables with $E|\xi_i|<\infty$. Let $E\xi_i=\mu $ and $S_n = \xi_1+\cdots +\xi_n$. Then
\begin{equation}
\displaystyle \frac{S_n}{n}\rightarrow \mu \ \mbox{a.s.} , n\rightarrow +\infty.
\end{equation}

\textbf{Theorem 6} Let $\xi_1,\xi_2,\cdots$ be i.i.d. random variables with $E\xi_i=0$ and $E\xi_i^2=\sigma^2<\infty$. Let $S_n = \xi_1 +\cdots+\xi_n$. If $\epsilon>0$, then
\begin{equation}
\displaystyle \frac{S_n}{\displaystyle n^{\frac{1}{2}}(\log n)^{\frac{1}{2}+\epsilon}}\rightarrow 0 \ \  \mbox{a.s.}
\end{equation}

\textbf{Theorem 7} (The De Moivre--Laplace CLT)
Let $\xi_1,\xi_2,\cdots$ be i.i.d. with $P(\xi_1=1)=P(\xi_1=-1)=\displaystyle  \frac{1}{2}$ and let $S_n=\xi_1+\cdots+\xi_n$. If $a<b$, then as $m\rightarrow \infty$,
\begin{equation}
\displaystyle  P(a \leq\frac{S_m}{\sqrt{m}}\leq b ) \rightarrow  \int_{a}^{b} \frac{1}{\sqrt{2\pi}}e^{-\frac{x^2}{2}}dx.
\end{equation}

\textbf{Theorem 8} (Lindeberg--L\'{e}vy  CLT) Suppose $\{\xi_1,\xi_2,\cdots\}$ is a sequence of i.i.d. random variables with $E(\xi_i) =\mu $ and $Var(\xi_i) = \sigma^2 < \infty$. Then as n approaches infinity, the random variables
$\displaystyle \frac{S_n-n\mu}{\sqrt{n}\sigma} $ converge in distribution to a normal $N(0,1)$:
\begin{equation}
\begin{array}{ll}
\displaystyle  \frac{S_{n}-n\mu}{\sqrt{n}\sigma}\ {\xrightarrow {d}}\ N(0,1).
\end{array}
\end{equation}

\textbf{Theorem 9} (Lindeberg CLT) Suppose $\{\xi_1, \xi_2,\cdots\}$ is a sequence of independent random variables, each with finite expected value $\mu_i$ and variance $\sigma^2_i$. Let
\begin{equation*}
\displaystyle B_{n}^{2}=\displaystyle \sum_{i=1}^{n}\sigma_{i}^{2}.
\end{equation*}
Suppose that for every $\epsilon>0$
\begin{equation}
{\displaystyle \lim_{n\rightarrow \infty }{\frac {1}{B_{n}^{2}}}\sum_{i=1}^{n}\operatorname {E} \left[(\xi_{i}-\mu_{i})^{2}\cdot \mathbf {1}_{\{|\xi_{i}-\mu_{i}|>\varepsilon B_{n}\}}\right]=0}
\end{equation}
where $\mathbf {1}_{\{\cdots\}}$ is the indicator function. Then the distribution of the standardized sums
\begin{equation}
\displaystyle \frac{1}{B_{n}} \sum_{i=1}^{n}\left(\xi_{i}-\mu_{i} \right) {\xrightarrow {d}}\ N(0,1).
\end{equation}

\textbf{Theorem 10} (Lyapunov CLT) Suppose $\{\xi_1, \xi_2,\cdots\}$ is a sequence of independent random variables, each with finite expected value $\mu_i$ and variance $\sigma^2_i$. Let
\begin{equation*}
\displaystyle B_{n}^{2}=\displaystyle \sum_{i=1}^{n}\sigma_{i}^{2}.
\end{equation*}
If for some $\delta>0$, Lyapunov¡¯s condition
\begin{equation}
\displaystyle \lim_{n\to \infty }{\frac{1}{B_{n}^{2+\delta}}}\sum_{i=1}^{n}\operatorname {E}\left[ |\xi_{i}- \mu_{i}|^{2+\delta } \right]=0
\end{equation}
is satisfied, then
\begin{equation}
\displaystyle \frac{1}{B_{n}}\sum_{i=1}^{n}\left(\xi_{i}-\mu_{i}\right)\ {\xrightarrow {d}}\ N(0,1).
\end{equation}

\textbf{Theorem 11} ( Lindeberg--Feller CLT) Suppose $\{\xi_1, \xi_2,\cdots\}$ is a sequence of independent random variables, $E(\xi_i)=\mu_i$ and $D(\xi_i)=\sigma^2_i$ . Let
\begin{equation}
\displaystyle B_{n}^{2}=\displaystyle \sum_{i=1}^{n}\sigma_{i}^{2}.
\end{equation}

Lindeberg condition
\begin{equation}
\displaystyle \lim_{n\rightarrow \infty} \frac{1}{B_{n}^{2}}\sum_{i=1}^{n}\operatorname {E}\left[|\xi_{i}- \mu_{i}|^{2} \cdot \mathbf {1}_{\{|\xi_{i}-\mu_{i}|>\varepsilon B_{n}\}} \right]=0
\end{equation}
holds iff
\begin{equation}
\displaystyle \frac{1}{B_{n}}\sum_{i=1}^{n}\left(\xi_{i}-\mu_{i}\right)\ {\xrightarrow {d}}\ N(0,1).
\end{equation}
and
\begin{equation}
 \displaystyle \lim\limits_{n\rightarrow \infty} \max\limits_{1\leq i\leq n} \frac{\sigma_{i}^{2}}{B_{n}^{2}}=0,
\end{equation}

For convenience, Denote
\begin{equation}
k^{*}=\left\{
\begin{array}{ll}
k,    & 0\leq k<n. \\
k(n), & 0\leq k(n)< n.\\
\end{array}
\right.
\  \ J_{k^{*}}=\left\{
\begin{array}{ll}
J_{k},   & k^{*}=k   \\
J_{k(n)},& k^{*}=k(n) \\
\end{array}
\right.
\end{equation}

If
\begin{equation}
\begin{array}{ll}
\displaystyle \lim \limits_{n\rightarrow \infty} \frac{k^{*}}{n}=0,
\end{array}
\end{equation}
we call it  asymptotic deleting negligibility condition of LLN.

If
\begin{equation}
\begin{array}{ll}
\displaystyle \lim \limits_{n\rightarrow \infty} \frac{k^{*}}{\sqrt{n}}=0,
\end{array}
\end{equation}
we call it  asymptotic deleting negligibility condition of CLT.

A simple example for $k(n)$ is $k(n)=[n^r]$, $0<r<1$.

\section{Deleting Items WLLN}
\label{}

\hspace*{1cm}

To demonstrate the deleting items WLLN can be proved by classical mathematical analysis methods, we discuss Bernoulli WLLN, Chebyshev WLLN, Khinchine WLLN and a general WLLN.

\subsection{Deleting Items Bernoulli WLLN}
\label{}

\textbf{Theorem 12} (Deleting Items Bernoulli WLLN) Let $\xi_{1}, \xi_{2}, \cdots, \xi_{n}, \cdots $ be mutual independent random sequence, where $E\xi_{i}=p$, $D\xi_{i}=p(1-p)$, $0<p<1$ . For $\forall \epsilon>0$, if natural number $k^*$ $(0<k^*<n)$ satisfies
\begin{equation*}
\begin{array}{ll}
\displaystyle  \lim\limits_{n \longrightarrow +\infty } \frac{k^*}{n}=0. \\
\end{array}
\end{equation*}
Then,
\begin{equation}
\begin{array}{ll}
\displaystyle  \lim\limits_{n \longrightarrow +\infty } P\left( \left|\frac{1}{n}\sum\limits_{i\in J\backslash J_{k^*}} \xi_i-p \right|< \varepsilon \right) = 1. \\
\end{array}
\end{equation}

\textbf{Proof}  Applying Chebyshev inequality, $\forall \epsilon>0$

\begin{equation*}
\begin{array}{ll}
\displaystyle P(| \frac{1}{n}\sum\limits_{i\in J\backslash J_{k^*}} \xi_i-\frac{1}{n}\sum\limits_{i\in J\backslash J_{k^*}}E\xi_i |\geq \epsilon) \leq \displaystyle \frac{D(\frac{1}{n}\sum\limits_{i\in J\backslash J_{k^*}} \xi_i )}{\epsilon^2} \\
\quad =\displaystyle \frac{(n-k^*) D(\xi_i)}{n^2 \epsilon^2}
 =\displaystyle \frac{(n-k^*) p(1-p)}{{n^2}\epsilon^2}\rightarrow 0, \ (n\rightarrow +\infty)
\end{array}
\end{equation*}

That is
\begin{equation*}
\begin{array}{ll}
\displaystyle \frac{1}{n}\sum\limits_{i\in J\backslash J_{k^*}} \xi_i-\frac{1}{n}\sum\limits_{i\in J\backslash J_{k^*}}E\xi_i \ {\xrightarrow {P}}\  0, \ \ n\rightarrow +\infty.
\end{array}
\end{equation*}

And,
\begin{equation*}
\begin{array}{ll}
\displaystyle  \frac{1}{n}\sum\limits_{i\in J_{k^*}} E\xi_i =\frac{k^*}{n} p \rightarrow 0, \ (n\rightarrow +\infty)
\end{array}
\end{equation*}

Then,
\begin{equation*}
\begin{array}{ll}
\displaystyle \frac{1}{n}\sum\limits_{i\in J\backslash J_{k^*}} \xi_i -p
=\frac{1}{n}\sum\limits_{i\in J\backslash J_{k^*}} \xi_i-\frac{1}{n}\sum\limits_{i\in J\backslash J_{k^*}}E\xi_i-\frac{1}{n}\sum\limits_{i\in J_{k^*}}E\xi_i \ {\xrightarrow {P}} \ 0, \ \ n\rightarrow +\infty.
\end{array}
\end{equation*}

which completes the proof.

\hfill $\Box$

\subsection{Deleting Items Chebyshev WLLN}
\label{}

Bernoulli WLLN is a special case of Chebyshev WLLN, the Deleting Items Chebyshev WLLN is given as follows.

\textbf{Theorem 13} (Deleting Items Chebyshev WLLN) Let $\xi_{1}, \xi_{2}, \cdots, \xi_{n}, \cdots $ be piecewise uncorrelated random sequence with bounded expectation and variance, where $ |E\xi_{i}|\leq M$, $ D\xi_{i}\leq C$, $i=1,2,\cdots $, $M$ and $C$ are constants.  For $\forall \epsilon>0$, if natural number $k^*$ $(0<k^*<n)$ satisfies
\begin{equation*}
\begin{array}{ll}
\displaystyle  \lim\limits_{n \longrightarrow +\infty } \frac{k^*}{n}=0. \\
\end{array}
\end{equation*}
Then,
\begin{equation}
\begin{array}{ll}
\displaystyle  \lim\limits_{n \rightarrow +\infty } P\left( \left|\frac{1}{n}\sum\limits_{i\in J\backslash J_{k^*}} \xi_i-\frac{1}{n} \sum\limits_{i=1}^n E\xi_i \right|< \varepsilon \right) = 1. \\
\end{array}
\end{equation}

\textbf{Proof}  According to Chebyshev inequality, $\forall \epsilon>0$

\begin{equation*}
\begin{array}{ll}
\displaystyle P(| \frac{1}{n}\sum\limits_{i\in J\backslash J_{k^*}} \xi_i-\frac{1}{n}\sum\limits_{i\in J\backslash J_{k^*}}E\xi_i |\geq \epsilon) \leq \displaystyle \frac{D(\frac{1}{n}\sum\limits_{i\in J\backslash J_{k^*}} \xi_i )}{\epsilon^2} \\
\quad =\displaystyle \frac{(n-k^*) D(\xi_i)}{n^2 \epsilon^2}
 \leq\displaystyle \frac{(n-k^*) C}{{n^2}\epsilon^2}\rightarrow 0, \ (n\rightarrow +\infty)
\end{array}
\end{equation*}

Hence,
\begin{equation*}
\begin{array}{ll}
\displaystyle \frac{1}{n}\sum\limits_{i\in J\backslash J_{k^*}} \xi_i-\frac{1}{n}\sum\limits_{i\in J\backslash J_{k^*}}E\xi_i \ {\xrightarrow {P}}\  0, \ \ n\rightarrow +\infty.
\end{array}
\end{equation*}
Again,
\begin{equation*}
\begin{array}{ll}
\displaystyle  |\frac{1}{n}\sum\limits_{i\in J_{k^*}} E\xi_i| \leq \frac{k^*}{n}M \rightarrow 0, \ (n\rightarrow +\infty)
\end{array}
\end{equation*}
Then,
\begin{equation*}
\begin{array}{ll}
\displaystyle  \frac{1}{n}\sum\limits_{i\in J_{k^*}} E\xi_i \rightarrow 0, \ (n\rightarrow +\infty)
\end{array}
\end{equation*}
Hence,
\begin{equation*}
\begin{array}{ll}
\displaystyle  \frac{1}{n} \sum\limits_{i\in J\backslash J_{k^*}} \xi_i -\frac{1}{n} \sum\limits_{i=1}^n E\xi_i
&= \displaystyle (\frac{1}{n} \sum\limits_{i\in J\backslash J_{k^*}} \xi_i-\frac{1}{n} \sum\limits_{i\in J\backslash J_{k^*}} E\xi_i)-\frac{1}{n} \sum\limits_{i\in J_{k^*}} E \xi_i \\
&\stackrel{P}{\longrightarrow} \ 0-0=0 \\
\end{array}
\end{equation*}

\hfill $\Box$

\subsection{Deleting Items Khinchine WLLN}
\label{}

\textbf{Theorem 14} (Deleting Items Khinchine WLLN) Let $\xi_{1}, \xi_{2}, \cdots, \xi_{n}, \cdots $ be independent identical distribution random variable sequence, with $E\xi_{i}=a$ $(i=1,2,\cdots)$. For $\forall \epsilon>0$, if natural number $k^*$ $(0<k^*<n)$ satisfies
\begin{equation*}
\begin{array}{ll}
\displaystyle  \lim\limits_{n \longrightarrow +\infty } \frac{k^*}{n}=0. \\
\end{array}
\end{equation*}
Then,
\begin{equation}
\begin{array}{ll}
\displaystyle  \frac{1}{n} \sum\limits_{i\in J\setminus J_{k^*}}\xi_i \stackrel{P}{\longrightarrow} a,  \  \  n\longrightarrow +\infty .\\
\end{array}
\end{equation}

\textbf{Proof} Since $\xi_{1}, \xi_{2}, \cdots, \xi_{n}, \cdots $ are independent identical distribution, they have same characteristic function $\varphi(t)$.
As $E\xi_i$ exists, then $\varphi(t)$ has expanded formula,
\begin{equation*}
\begin{array}{ll}
\displaystyle  \varphi(t) = \varphi(0)+\varphi'(0) t + o(t)=1+iat+o(t)
\end{array}
\end{equation*}
According to the independent property, the characteristic function of $\displaystyle \frac{1}{n}\sum\limits_{i\in J\setminus J_{k^*}} \xi_i$ is
\begin{equation*}
\begin{array}{ll}
\displaystyle  [\varphi(\frac{t}{n})]^{n-{k^*}} =\displaystyle [1+ia\frac{t}{n}+o(\frac{t}{n})]^{n-{k^*}}
\end{array}
\end{equation*}
For any $t$ ,
\begin{equation*}
\begin{array}{ll}
\displaystyle \lim\limits_{n \longrightarrow +\infty } [\varphi(\frac{t}{n})]^{n-{k^*}} =\displaystyle \lim\limits_{n \longrightarrow +\infty }[1+ia\frac{t}{n}+o(\frac{t}{n})]^{n\frac{n-{k^*}}{n}} = \displaystyle e^{iat}
\end{array}
\end{equation*}
As $e^{iat}$ is the characteristic function of distribution function
$$F(x)=\displaystyle \left\{ \begin{array}{cc} 1,& x>a \\0,& x\leq a \end{array} \right. $$
And, it is the distribution function of random variable $\eta =a$.
Hence, the distribution function of $ \frac{1}{n}\sum\limits_{i\in J\setminus J_k} \xi_i $ is weakly convergent to $F(x)$, and
\begin{equation*}
\displaystyle \frac{1}{n}\sum\limits_{i\in J\setminus J_{k^*}} \xi_i \stackrel{P}{\longrightarrow} a.
\end{equation*}

\hfill $\Box$

In fact, an alternative simple proof of Deleting Items Khinchine WLLN is as follows.

\textbf{Proof} According to Theorem 3 Khinchine WLLN,
\begin{equation*}
\displaystyle \frac{1}{n}\sum\limits_{i\in J\setminus J_{k^*}} \xi_i=\displaystyle \frac{n-k^*}{n} \frac{1}{n-k^*}\sum\limits_{i\in J\setminus J_{k^*}} \xi_i \stackrel{P}{\longrightarrow} 1\cdot a =a.
\end{equation*}

\hfill $\Box$

\subsection{General Deleting items WLLN}
\label{}

According to the theorem of convergence in probability and convergence in squared mean, we have a general deleting item theorem.

\textbf{Theorem 15} (General Deleting items WLLN) Let $\xi_{1}, \xi_{2}, \cdots, \xi_{n}, \cdots $ be a random variable sequence with bounded expectation $|E\xi_{i}|<M $, $i=1,2,\cdots $, where $M$ is a constant. For $\forall \epsilon>0$, if natural number $k^*$ $(0<k^*<n)$ satisfies
\begin{equation*}
\begin{array}{ll}
\displaystyle  \lim\limits_{n \longrightarrow +\infty } \frac{k^*}{n}=0. \\
\end{array}
\end{equation*}
Then,
\begin{equation}
\begin{array}{ll}
\displaystyle \frac{1}{n} \sum\limits_{i\in J\setminus J_{k^*}} \xi_i-\frac{1}{n} \sum\limits_{i=1}^n E\xi_i  \stackrel{P}{\longrightarrow} 0 , \ \ n \longrightarrow +\infty.\\
\Longleftrightarrow \displaystyle E[\frac{|\frac{1}{n} \sum\limits_{i\in J\setminus J_{k^*}} \xi_i-\frac{1}{n} \sum\limits_{i=1}^n E\xi_i|^2}
{1+|\frac{1}{n} \sum\limits_{J\setminus J_{k^*}} \xi_i-\frac{1}{n} \sum\limits_{i=1}^n E\xi_i|^2}  ] \stackrel{P}{\longrightarrow} 0 , \ \ n \longrightarrow +\infty.
\end{array}
\end{equation}

\textbf{Proof}
Since $|E\xi_{i}|<M $, $i=1,2,\cdots $, we obtain
\begin{equation*}
\begin{array}{ll}
\displaystyle  |\frac{1}{n}\sum\limits_{i\in J_{k^*}} E\xi_i| \leq \frac{k^*}{n}M \rightarrow 0, \ (n\rightarrow +\infty).
\end{array}
\end{equation*}
That is,
\begin{equation*}
\begin{array}{ll}
\displaystyle  \frac{1}{n}\sum\limits_{i\in J_{k^*}} E\xi_i \rightarrow 0, \ (n\rightarrow +\infty).
\end{array}
\end{equation*}

If
\begin{equation*}
\begin{array}{ll}
\displaystyle \frac{1}{n} \sum\limits_{i\in J} \xi_i - \frac{1}{n}\sum\limits_{i\in J}E\xi_i \stackrel{P}{\longrightarrow} 0 \\
\end{array}
\end{equation*}
then
\begin{equation*}
\begin{array}{ll}
\displaystyle \frac{1}{n} \sum\limits_{i\in J\setminus J_{k^*}}\xi_i - \frac{1}{n}\sum\limits_{i\in J} E\xi_i \\
\displaystyle =\frac{n-k^*}{n} [\frac{1}{n-k^*} \sum\limits_{i\in J\setminus J_{k^*}} \xi_i - \frac{1}{n-k^*} \sum\limits_{i\in J\setminus J_{k^*}} E\xi_i]-\frac{1}{n} \sum\limits_{i\in J_{k^*}} E\xi_i\\
\stackrel{P}{\longrightarrow} 1\cdot 0+0=0 \\
\end{array}
\end{equation*}
According to Theorem 4,
\begin{equation*}
\begin{array}{ll}
\displaystyle \frac{1}{n} \sum\limits_{i\in J\setminus J_{k^*}} \xi_i - \frac{1}{n} \sum\limits_{i\in J\setminus J_{k^*}} E\xi_i-\frac{1}{n} \sum\limits_{i\in J_{k^*}} E\xi_i \stackrel{P}{\longrightarrow} 0 \\
\Longleftrightarrow  \displaystyle E[\frac{|\frac{1}{n} \sum\limits_{i\in J\setminus J_{k^*}} \xi_i-\frac{1}{n} \sum\limits_{i\in J\setminus J_{k^*}} E\xi_i- \frac{1}{n} \sum\limits_{i\in J_{k^*}} E\xi_i|^2}
{1+|\frac{1}{n} \sum\limits_{i\in J\setminus J_{k^*}} \xi_i-\frac{1}{n} \sum\limits_{i\in J\setminus J_{k^*}} E\xi_i- \frac{1}{n} \sum\limits_{i\in J_{k^*}} E\xi_i|^2}] \stackrel{P}{\longrightarrow} 0 \\
\Longleftrightarrow \displaystyle E |\frac{|\frac{1}{n} \sum\limits_{i\in J\setminus J_{k^*}} \xi_i-\frac{1}{n} \sum\limits_{i=1}^n E\xi_i|^2}
{1+|\frac{1}{n}\sum\limits_{J\setminus J_{k^*}} \xi_i-\frac{1}{n} \sum\limits_{i=1}^n E\xi_i|^2}] \stackrel{P}{\longrightarrow} 0 \\
\end{array}
\end{equation*}

Hence, the formula(32) holds.

\hfill $\Box$

\section{Uniform of Deleting items WLLN}

In the above section, we strictly examine the proof of the classical WLLN in literatures, a uniform Deleting items WLLN could be given by Slutsky' Theorem for all three WLLN as follows.

\textbf{Theorem 16} Let $\xi_{1}, \xi_{2}, \cdots, \xi_{n}, \cdots $ be a random variable sequence with bounded expectation $|E\xi_{i}|<M $, $i=1,2,\cdots $, where $M$ is a constant. For natural number $k^*$ $(0<k^*<n)$ satisfies
\begin{equation*}
\begin{array}{ll}
\displaystyle  \lim\limits_{n \longrightarrow +\infty } \frac{k^*}{n}=0. \\
\end{array}
\end{equation*}
If
\begin{equation*}
\begin{array}{ll}
\displaystyle \frac{1}{n} \sum\limits_{i\in J} \xi_i-\frac{1}{n} \sum\limits_{i=1}^n E\xi_i  \stackrel{P}{\longrightarrow} 0 , \ \ n \longrightarrow +\infty.\\
\end{array}
\end{equation*}
Then,
\begin{equation}
\begin{array}{ll}
\displaystyle \frac{1}{n} \sum\limits_{i\in J\setminus J_{k^*}} \xi_i-\frac{1}{n} \sum\limits_{i=1}^n E\xi_i  \stackrel{P}{\longrightarrow} 0 , \ \ n \longrightarrow +\infty.\\
\end{array}
\end{equation}

\textbf{Proof}
Since
\begin{equation*}
\begin{array}{ll}
\displaystyle \frac{1}{n} \sum\limits_{i\in J\setminus J_{k^*}} \xi_i-\frac{1}{n} \sum\limits_{i=1}^n E\xi_i  \\
=\displaystyle \frac{1}{n} \sum\limits_{i\in J\setminus J_{k^*}} \xi_i-\frac{1}{n}\sum\limits_{i\in J\setminus J_{k^*}} E\xi_i- \frac{1}{n}\sum\limits_{i\in J_{k^*}} E\xi_i  \\
=\displaystyle \frac{n-k^*}{n} (\frac{1}{n-k^*} \sum\limits_{i\in J\setminus J_{k^*}} \xi_i-\frac{1}{n-k^*}\sum\limits_{i\in J\setminus J_{k^*}} E\xi_i)- \frac{1}{n}\sum\limits_{i\in J_{k^*}} E\xi_i  \\
\stackrel{P}{\longrightarrow} 1\cdot 0+0=0 , \ \ n \longrightarrow +\infty.\\
\end{array}
\end{equation*}
it ends the proof.

\hfill $\Box$

\section{Deleting items SLLN}
\label{}

In this section ,we give two deleting items SLLN theorems.

\textbf{Theorem 17} (Deleting items SLLN). Let $\xi_1,\xi_2,\cdots$ be pairwise independent identically distributed ( i.i.d.) random variables with $E|\xi_i|<\infty$. Let $E\xi_i=\mu $ and $S_n=\xi_1+\cdots+\xi_n$. If natural number $k^*$ $(0<k^*<n)$ satisfies
\begin{equation*}
\begin{array}{ll}
\displaystyle  \lim\limits_{n \longrightarrow +\infty } \frac{k^*}{n}=0. \\
\end{array}
\end{equation*}
Then,
\begin{equation}
\displaystyle \frac{S_{J\setminus J_{k^{*}}}}{n}\rightarrow \mu \ \mbox{a.s.} ,\  n\rightarrow +\infty.
\end{equation}

\textbf{Proof} Applying Theorem 5, we obtain
\begin{equation*}
\begin{array}{ll}
\displaystyle \frac{S_{J\backslash J_{k^{*}}}}{n}
&=\displaystyle \frac{S_{J\backslash J_{k^{*}}}}{(n-k^{*})} {\frac{(n-k^{*})}{n}} \\
&\displaystyle {\xrightarrow {a.s}}\ \mu \cdot 1= \mu,
\end{array}
\end{equation*}
Hence, it ends the proof.

\hfill $\Box$

\textbf{Theorem 18} Let $\xi_1,\xi_2,\cdots$ be i.i.d. random variables with $E\xi_i=0$ and $E\xi_i^2=\sigma^2<\infty$. Let $S_n = \xi_1+\cdots+\xi_n$. If $\epsilon>0$, and natural number $k^*$ $(0<k^*<n)$ satisfies
\begin{equation*}
\begin{array}{ll}
\displaystyle  \lim\limits_{n \longrightarrow +\infty } \frac{k^*}{n}=0, \\
\end{array}
\end{equation*}
then
\begin{equation}
\begin{array}{ll}
\displaystyle \frac{S_{J\setminus J_{k^*}}}{\displaystyle n^{\frac{1}{2}}(\log n)^{\frac{1}{2}+\epsilon}}\rightarrow 0 \ \  \mbox{a.s.} \\
\end{array}
\end{equation}

\textbf{Proof} Since
\begin{equation*}
\begin{array}{ll}
\displaystyle \lim \limits_{n\rightarrow \infty} \frac{n-k^*}{n}=\lim \limits_{n\rightarrow \infty} 1- \frac{k^{*}}{n}=1 \\
\displaystyle \lim \limits_{n\rightarrow \infty} \frac{\log(n-k^*)}{\log n}=\lim\limits_{n\rightarrow \infty} \frac{\log n + \log(1-\frac{k^{*}}{n})}{\log n}=\lim \limits_{n\rightarrow \infty} 1+ \frac{ \log(1-\frac{k^{*}}{n})}{\log n}=1
\end{array}
\end{equation*}
and according to Theorem 6, we obtain
\begin{equation*}
\begin{array}{ll}
\displaystyle \frac{S_{J\backslash J_{k^{*}}} } {\displaystyle n^{\frac{1}{2}}(\log n)^{\frac{1}{2}+\epsilon}}
&=\displaystyle \frac{S_{J\backslash J_{k^{*}}}}{(n-k^{*})^{\frac{1}{2}}(\log (n-k^{*}))^{\frac{1}{2}+\epsilon}}  \frac{(n-k^{*})^{\frac{1}{2}}}{n^{\frac{1}{2}}}\frac{(\log (n-k^{*}))^{\frac{1}{2}+\epsilon}}{(\log n)^{\frac{1}{2}+\epsilon}} \\
&\displaystyle {\xrightarrow {a.s}}\ 0\cdot1\cdot1 = 0.
\end{array}
\end{equation*}
Hence, it completes the proof.

\hfill $\Box$

\section{Deleting items CLT}
\label{}

Since De Moivre-Laplace CLT is a special case of Lindeberg-L\'{e}vy CLT, Lyapunov CLT is a special case of Lindeberg CLT (a sequence satisfies Lyapunov condition, it satisfies Lindeberg condition), and Lindeberg-Feller CLT is a necessary and sufficient theorem, it is stronger than Lyapunov CLT, we only prove deleting items Lindeberg--L\'{e}vy CLT and deleting items Lindeberg--Feller CLT.

\textbf{Theorem 19} (Deleting items De Moivre--Laplace CLT)
Let $\xi_1,\xi_2,\cdots$ be i.i.d. with $P(\xi_1=1)=P(\xi_1=-1)=\displaystyle \frac{1}{2}$ and let $S_n=\xi_1+\cdots+\xi_n$. If $a<b$, and natural number $k^*$ $(0<k^*<n)$ satisfies
\begin{equation*}
\begin{array}{ll}
\displaystyle  \lim\limits_{n \longrightarrow +\infty } \frac{k^*}{n}=0, \\
\end{array}
\end{equation*}
then
\begin{equation}
\begin{array}{ll}
\displaystyle  P(a \leq\frac{S_{J\setminus J_{k^*}}}{\sqrt{n}}\leq b )\rightarrow \displaystyle  \int_{a}^{b} \frac{1}{\sqrt{2\pi}}e^{-\frac{x^2}{2}}dx, \ \ n\rightarrow \infty. \\
\end{array}
\end{equation}

\textbf{Theorem 20} (Deleting items Lindeberg--L\'{e}vy CLT) Suppose $\{\xi_1,\xi_2,\cdots\}$ is a sequence of i.i.d. random variables with $E(\xi_i) =\mu $ and $D(\xi_i) = \sigma^2 < \infty$. If natural number $k^*$ $(0<k^*<n)$ satisfies
\begin{equation*}
\begin{array}{ll}
\displaystyle \lim \limits_{n\rightarrow \infty} \frac{k^{*}}{\sqrt{n}}=0
\end{array}
\end{equation*}
then,
\begin{equation}
\begin{array}{ll}
\displaystyle  \frac{S_{J\setminus J_{k^*}}-n\mu}{\sqrt{n}\sigma}\ {\xrightarrow {d}}\ N(0,1). \\
\end{array}
\end{equation}

\textbf{Proof}: Since
\begin{equation*}
\begin{array}{ll}
\displaystyle  \frac{S_{n}-n\mu}{\sqrt{n}\sigma} \ {\xrightarrow {d}}\ N\left(0,1\right),
\end{array}
\end{equation*}
we obtain
\begin{equation*}
\begin{array}{ll}
\displaystyle \frac{\sum\limits_{J\backslash J_{k^{*}}} \xi_{i}-n\mu } {\sqrt{n}\sigma}
=\displaystyle \frac{\sum\limits_{J\backslash J_{k^{*}}} \xi_{i}-(n-k^{*})\mu-k^{*}\mu } {\sqrt{n}\sigma}\\
=\displaystyle \frac{\sum\limits_{J\backslash J_{k^{*}}} \xi_{i}-(n-k^{*})\mu} {\displaystyle \sqrt{n-k^{*}}\sigma} \frac{\displaystyle \sqrt{n-k^{*}}}{\sqrt{n}}-\frac{k^{*}\mu}{\sqrt{n}\sigma}\\
\displaystyle {\xrightarrow {d}}\ N(0,1)\cdot 1+0 {\xrightarrow {d}}\ N(0,1).
\end{array}
\end{equation*}
Hence, it ends the proof.

\hfill $\Box$

From the proof, when $\mu=0$, the convergent condition satisfying $\displaystyle  \lim\limits_{n \longrightarrow +\infty } \frac{k^*}{n}=0$  is enough. And Deleting items De Moivre--Laplace CLT takes the special case with $\mu=0$. When $\mu\neq 0$, the convergent condition of $k^*$ is stronger than that in case of $\mu=0$.

\textbf{Theorem 21} (Deleting items Lindeberg CLT) Suppose $\{\xi_1, \xi_2,\cdots\}$ is a sequence of independent random variables, each with finite expected value $\mu_i$ and variance $\sigma^2_i$. Let
\begin{equation*}
\displaystyle B_{n}^{2}=\displaystyle \sum_{i=1}^{n}\sigma_{i}^{2}.
\end{equation*}
Suppose that for every $\epsilon>0$
\begin{equation*}
{\displaystyle \lim_{n\rightarrow \infty }{\frac {1}{B_{n}^{2}}}\sum_{i=1}^{n}\operatorname {E} \left[(\xi_{i}-\mu_{i})^{2}\cdot \mathbf {1}_{\{|\xi_{i}-\mu_{i}|>\varepsilon B_{n}\}}\right]=0}
\end{equation*}
where $\mathbf {1}_{\{\cdots\}}$ is the indicator function.
If natural number $k^*$ $(0<k^*<n)$ satisfies
\begin{equation*}
\begin{array}{ll}
\displaystyle \lim \limits_{n\rightarrow \infty} \frac{k^{*}}{\sqrt{n}}=0
\end{array}
\end{equation*}
and,
\begin{equation*}
\begin{array}{ll}
\displaystyle \frac{ \max \limits_{1\leq i\leq n} \sigma_i^2 }{B_{n}^2}=O(\frac{1}{n})\\
\displaystyle \frac{ \max \limits_{1\leq i\leq n} |\mu_i| }{B_{n}}=O(\frac{1}{\sqrt{n}})
\end{array}
\end{equation*}
Then the distribution of the standardized sums
\begin{equation}
\displaystyle \frac{1}{B_{n}} \left( S_{J\backslash J_{k^{*}}}-\sum_{i=1}^n \mu_{i} \right){\xrightarrow {d}} N(0,1)
\end{equation}

\textbf{Theorem 22} (Deleting items Lyapunov CLT) Suppose $\{\xi_1, \xi_2,\cdots\}$ is a sequence of independent random variables, each with finite expected value $E\xi_i=\mu_i$ and variance $D\xi_i=\sigma^2_i$. Define
\begin{equation*}
\displaystyle B_{n}^{2}=\displaystyle \sum_{i=1}^{n}\sigma_{i}^{2}
\end{equation*}
If for some $\delta>0$, Lyapunov¡¯s condition
\begin{equation*}
\displaystyle \lim_{n\to \infty }{\frac{1}{B_{n}^{2+\delta}}}\sum_{i=1}^{n}\operatorname {E}\left[ |\xi_{i}- \mu_{i}|^{2+\delta } \right]=0
\end{equation*}
is satisfied.
If natural number $k^*$ $(0<k^*<n)$ satisfies
\begin{equation*}
\begin{array}{ll}
\displaystyle \lim \limits_{n\rightarrow \infty} \frac{k^{*}}{\sqrt{n}}=0
\end{array}
\end{equation*}
and,
\begin{equation*}
\begin{array}{ll}
\displaystyle \frac{ \max \limits_{1\leq i\leq n} \sigma_i^2 }{B_{n}^2}=O(\frac{1}{n})\\
\displaystyle \frac{ \max \limits_{1\leq i\leq n} |\mu_i| }{B_{n}}=O(\frac{1}{\sqrt{n}})
\end{array}
\end{equation*}
then
\begin{equation}
\displaystyle \frac{1}{B_{n}} \left( S_{J\backslash J_{k^{*}}}-\sum_{i=1}^n \mu_{i} \right){\xrightarrow {d}} N(0,1)
\end{equation}

\textbf{Theorem 23} (Deleting Items Lindeberg-Feller CLT)
Suppose $\{\xi_1, \xi_2,\cdots\}$ is a sequence of independent random variables, $E(\xi_i)=\mu_i <+\infty$ and $Var(\xi_i)=\sigma^2_i <+\infty $ . Define
\begin{equation}
\displaystyle B_{n}^{2}=\displaystyle \sum_{i=1}^{n}\sigma_{i}^{2}
\end{equation}
If the following conditions hold,

(1) Lindeberg condition
\begin{equation}
\displaystyle \lim_{n\rightarrow \infty} \frac{1}{B_{n}^{2}}\sum_{i=1}^{n}\operatorname {E}\left[|\xi_{i}- \mu_{i}|^{2} \cdot \mathbf {1}_{\{|\xi_{i}-\mu_{i}|>\varepsilon B_{n}\}} \right]=0
\end{equation}

(2)
\begin{equation}
\begin{array}{ll}
\displaystyle \frac{ \max \limits_{1\leq i\leq n} \sigma_i^2 }{B_{n}^2}=O(\frac{1}{n})
\end{array}
\end{equation}

(3)
\begin{equation}
\begin{array}{ll}
\displaystyle \frac{ \max \limits_{1\leq i\leq n} |\mu_i| }{B_{n}}=O(\frac{1}{\sqrt{n}})
\end{array}
\end{equation}

(4) For natural number $k^*$ $(0<k^*<n)$ satisfies
\begin{equation*}
\begin{array}{ll}
\displaystyle \lim \limits_{n\rightarrow \infty} \frac{k^{*}}{\sqrt{n}}=0
\end{array}
\end{equation*}
then
\begin{equation}
\displaystyle \frac{1}{B_{n}}\left(S_{J\backslash J_{k^{*}}}- \sum_{i=1}^{n}\mu_{i}\right)\ {\xrightarrow {d}}\ N(0,1).
\end{equation}

\textbf{Proof}
If
\begin{equation*}
\begin{array}{ll}
\displaystyle \lim \limits_{n\rightarrow \infty} \frac{k^{*}}{\sqrt{n}}=0,
\end{array}
\end{equation*}
then
\begin{equation*}
\begin{array}{ll}
\displaystyle \lim \limits_{n\rightarrow \infty} \frac{k^{*}}{n}=0.
\end{array}
\end{equation*}

Since
\begin{equation*}
\displaystyle \lim_{n\rightarrow \infty} \frac{1}{B_{n}^{2}}\sum_{i=1}^{n}\operatorname {E}\left[|\xi_{i}- \mu_{i}|^{2} \cdot \mathbf {1}_{\{|\xi_{i}-\mu_{i}|>\varepsilon B_{n}\}} \right]=0
\end{equation*}
According to Lindeberg--Feller CLT, we obtain
\begin{equation*}
\displaystyle \frac{1}{B_{J\backslash J_{k^{*}}}}\left(S_{J\backslash J_{k^{*}}}- \sum_{J\backslash J_{k^{*}}}\mu_{i}\right) {\xrightarrow {d}}\ N(0,1).
\end{equation*}
Since
\begin{equation*}
\begin{array}{ll}
\displaystyle | \frac{B_{J_{k^{*}}}^2}{B_{J}^2} |
\leq \displaystyle  \frac{k^{*}\max \limits_{1\leq i\leq n} \sigma_i^2 }{B_{J}^2}
= \displaystyle  \frac{k^{*}}{n}  n  \frac{\max \limits_{1\leq i\leq n} \sigma_i^2 }{B_{J}^2}
= \displaystyle  \frac{k^{*}}{n}   ( \frac{\max \limits_{1\leq i\leq n} \sigma_i^2 }{B_{J}^2}/\frac{1}{n})
\longrightarrow 0 , n\longrightarrow +\infty. \\
\end{array}
\end{equation*}
we obtain
\begin{equation*}
\begin{array}{ll}
\displaystyle \frac{B_{J\backslash J_{k^{*}}}}{B_{J}}=\displaystyle 1-\frac{B_{J_{k^{*}}}}{B_{J}}\longrightarrow 1 , \ \  n\longrightarrow +\infty. \\
\end{array}
\end{equation*}
And
\begin{equation*}
\begin{array}{ll}
\displaystyle | \frac{1}{B_{n}} \sum_{J_{k^{*}}}\mu_{i} |
&\leq \displaystyle  \frac{k^{*}\max \limits_{1\leq i\leq n} |\mu_i |}{B_{J}}
= \displaystyle  \frac{k^{*}}{\sqrt{n}}  \sqrt{n}  \frac{\max \limits_{1\leq i\leq n} | \mu_i| }{B_{J}}\\
&= \displaystyle  \frac{k^{*}}{\sqrt{n}} (\frac{\max \limits_{1\leq i\leq n} | \mu_i| }{B_{J}}/\frac{1}{\sqrt{n}})
\longrightarrow 0 ,\  n\longrightarrow +\infty. \\
\end{array}
\end{equation*}
Then,
\begin{equation*}
\begin{array}{ll}
\displaystyle \frac{1}{B_{n}}\left(S_{J\backslash J_{k^{*}}}- \sum_{i=1}^{n}\mu_{i}\right)
=\displaystyle \frac{1}{B_{n}}\left(S_{J\backslash J_{k^{*}}}- \sum_{J\backslash J_{k^{*}}}\mu_{i}-\sum_{J_{k^{*}}}\mu_{i}\right) \\
=\displaystyle \frac{B_{J\backslash J_{k^{*}}}}{B_{J}}\frac{1}{B_{J\backslash J_{k^{*}}}}\left(S_{J\backslash J_{k^{*}}}- \sum_{J\backslash J_{k^{*}}}\mu_{i}\right) -\frac{1}{B_{n}} \sum_{J_{k^{*}}}\mu_{i} \\
=\displaystyle \frac{B_{J\backslash J_{k^{*}}}}{B_{J}}\frac{1}{B_{J\backslash J_{k^{*}}}}\left(\sum_{J\backslash J_{k^{*}}} (\xi_i-\mu_{i})\right) -\frac{1}{B_{n}} \sum_{J_{k^{*}}}\mu_{i} \\
{\xrightarrow {d}}\ N(0,1).
\end{array}
\end{equation*}

Hence, it ends the proof.

\hfill $\Box$

Note that if $\mu_1=\mu_2=\cdots=\mu_n$, $\sigma_1=\sigma_2=\cdots=\sigma_n$, then
\begin{equation*}
\begin{array}{ll}
\displaystyle \frac{B_{J\backslash J_{k^{*}}}}{B_{J}}=\displaystyle \frac{n-k^{*}}{n}, \\
\displaystyle \frac{1}{B_{n}} \sum_{J_{k^{*}}}\mu_{i}=\displaystyle \frac{k^{*}}{\sqrt{n}}\frac{\mu}{\sigma}.\\
\end{array}
\end{equation*}
Therefore, it is easy to understand the different convergent condition of $k^{*}$ in each deleting items CLT.

\section{Application of asymptotic estimate of deleting items WLLN}
\label{}

The deleting items WLLN theorems can explain the phenomena why we can get almost same convergent results as the classical WLLN conclusions, though we do not strictly constrict the experiments condition in practical performance.

Obviously, $\displaystyle \frac{1}{n}\sum\limits_{i\in J\backslash J_{k^*}} \xi_i$ is not the subsequence of  $\displaystyle \frac{1}{n}\sum\limits_{i\in J} \xi_i$. And, generally
$$\displaystyle P(\frac{1}{n}\sum\limits_{i\in J\backslash J_{k^*}} \xi_i \neq \frac{1}{n}\sum\limits_{i\in J} \xi_i)=1,$$\\
especially when $\xi_i$ ($i=1,2,\cdots$) are continuous random variables.

In fact,
\begin{equation*}
\begin{array}{ll}
\displaystyle P(\frac{1}{n}\sum\limits_{i\in J\backslash J_{k^*}} \xi_i \neq \frac{1}{n}\sum\limits_{i\in J} \xi_i)
=\displaystyle P(\frac{1}{n}\sum\limits_{i\in J} \xi_i-\frac{1}{n}\sum\limits_{i\in J\backslash J_{k^*}} \xi_i \neq 0)\\
=\displaystyle P(\frac{1}{n}\sum\limits_{i\in J_{k^*}} \xi_i \neq 0)=1
\end{array}
\end{equation*}

Theoretically, WLLN is the foundation of moment estimator, bias and consistency are important concepts relative to estimator in sampling theory of statistics [17-18].
Suppose $\xi_{1}, \xi_{2}, \cdots, \xi_{n}, \cdots $ are i.i.d random  variables of  $\xi$ with $E\xi=\mu$ , $D\xi=\sigma^2$. Denote
\begin{equation}
\begin{array}{ll}
\displaystyle \overline{X}=\displaystyle  \frac{1}{n}\sum_{i=1}^n \xi_i , \ \
\displaystyle S^2= \frac{1}{n} \sum_{i=1}^n (\xi_i-\overline{X})^2.
\end{array}
\end{equation}
Then,
\begin{equation}
\begin{array}{ll}
\displaystyle E\overline{X}=\mu , \ \  ES^2= \sigma^2.
\end{array}
\end{equation}

The consistency of $\displaystyle \frac{1}{n} \sum\limits_{i\in J\setminus J_{k}} (\xi_i-\bar{\xi})^2$ and $\displaystyle \frac{1}{n} \sum\limits_{i\in J\setminus J_{k(n)}} (\xi_i-\bar{\xi})^2$ are discussed in deleting item WLLN. The  discussion of bias is in the following theorem.

\textbf{Theorem 24} Suppose $\xi_{1}, \xi_{2}, \cdots, \xi_{n}, \cdots $ are i.i.d random  variables of  $\xi$,  with $E\xi=\mu$, $D\xi=\sigma^2$. Denote
\begin{equation}
\begin{array}{ll}
\displaystyle \widetilde{X}=\displaystyle \frac{1}{n} \sum\limits_{i\in J\setminus J_{k^*}} \xi_i, &
\displaystyle \widetilde{S}_{1}^2=\displaystyle \frac{1}{n}\sum\limits_{i=1}^n (\xi_i-\widetilde{X})^2,\\
\displaystyle \widetilde{S}_{2}^2=\displaystyle \frac{1}{n}\sum\limits_{i\in J\setminus J_{k^*}} (\xi_i-\overline{X})^2,&
\displaystyle \widetilde{S}_{3}^2=\displaystyle \frac{1}{n}\sum\limits_{i\in J\setminus J_{k^*}} (\xi_i-\widetilde{X})^2.\\
\end{array}
\end{equation}
Then,
\begin{equation}
\begin{array}{ll}
\displaystyle \widetilde{S}_{1}^2=\displaystyle \frac{1}{n}\sum\limits_{i=1}^n \xi_{i}^2-2\widetilde{X}\overline{X}+\widetilde{X}^2,\\
\displaystyle \widetilde{S}_{2}^2=\displaystyle \frac{1}{n}\sum\limits_{i\in J\setminus J_{k^*}} \xi_{i}^2-2\widetilde{X}\overline{X}+(1-\frac{k^{*}}{n})\overline{X}^2,\\
\displaystyle \widetilde{S}_{3}^2=\displaystyle \frac{1}{n}\sum\limits_{i\in J\setminus J_{k^*}} \xi_i^2-(1+\frac{k^{*}}{n})\widetilde{X}^2.\\
\displaystyle E\widetilde{X}=\displaystyle (1-\frac{k^*}{n})\mu , \\
\displaystyle E\widetilde{S}_{1}^2=\displaystyle (1-\frac{1}{n}+\frac{k^{*}}{n^2})\sigma^2+\frac{{k^*}^2}{n^2}\mu^2.\\
\displaystyle E\widetilde{S}_{2}^2=\displaystyle (1-\frac{1}{n}-\frac{k^{*}}{n}+\frac{k^{*}}{n^2})\sigma^2\\
\displaystyle E\widetilde{S}_{3}^2=\displaystyle (1-\frac{1}{n}-\frac{k^{*}}{n}+\frac{k^{*}}{n^3})\sigma^2+(1-\frac{k^*}{n})\frac{k^*}{n^2}\mu^2\\
\end{array}
\end{equation}

\textbf{Proof}

(1) Simply expanding the sum, we obtain
\begin{equation*}
\begin{array}{ll}
\displaystyle \widetilde{S}_{1}^2&=\displaystyle \frac{1}{n}\sum\limits_{i=1}^n (\xi_i-\widetilde{X})^2
=\displaystyle \frac{1}{n}\sum\limits_{i=1}^n (\xi_i^2-2 \widetilde{X} \xi_i+\widetilde{X}^2)\\
&=\displaystyle \frac{1}{n}\sum\limits_{i=1}^n \xi_{i}^2-2\widetilde{X}\overline{X}+\widetilde{X}^2,\\
\end{array}
\end{equation*}

\begin{equation*}
\begin{array}{ll}
\displaystyle \widetilde{S}_{2}^2&=\displaystyle \frac{1}{n}\sum\limits_{i\in J\setminus J_{k^*}} (\xi_i-\overline{X})^2=\displaystyle \frac{1}{n}\sum\limits_{i\in J\setminus J_{k^*}} (\xi_i^2-2\overline{X}\xi_i +\overline{X}^2)\\
&=\displaystyle \frac{1}{n}\sum\limits_{i\in J\setminus J_{k^*}} \xi_{i}^2-2\widetilde{X}\overline{X}+(1-\frac{k^{*}}{n})\overline{X}^2,\\
\end{array}
\end{equation*}

\begin{equation*}
\begin{array}{ll}
\displaystyle \widetilde{S}_{3}^2 &=\displaystyle \frac{1}{n}\sum\limits_{i\in J\setminus J_{k^*}} (\xi_i-\widetilde{X})^2=\displaystyle \frac{1}{n}\sum\limits_{i\in J\setminus J_{k^*}} (\xi_i^2-2\widetilde{X}\xi_i+\widetilde{X}^2)\\
&=\displaystyle \frac{1}{n}\sum\limits_{i\in J\setminus J_{k^*}} \xi_i^2-(1+\frac{k^{*}}{n})\widetilde{X}^2.\\
\end{array}
\end{equation*}

(2) Since
\begin{equation*}
\begin{array}{ll}
\displaystyle E\widetilde{X}&=\displaystyle E(\frac{1}{n} \sum\limits_{i\in J\setminus J_{k^*}} \xi_i)=\displaystyle (1-\frac{k^*}{n})\mu ,
\end{array}
\end{equation*}

\begin{equation*}
\begin{array}{ll}
\displaystyle E(\widetilde{X}\overline{X})&=\displaystyle E((\frac{1}{n}\sum\limits_{i\in J\setminus J_{k^*}} \xi_i ) (\frac{1}{n} \sum\limits_{i=1}^n \xi_i))\\
&=\displaystyle \frac{1}{n^2}E( \sum\limits_{i\in J\setminus J_{k^*}} (\xi_i+\mu-\mu) \sum\limits_{i=1}^n (\xi_i+\mu-\mu))\\
&=\displaystyle (\frac{1}{n}-\frac{k^*}{n^2})\sigma^2+(1-\frac{k^*}{n})\mu^2, \\
\end{array}
\end{equation*}

\begin{equation*}
\begin{array}{ll}
\displaystyle E(\overline{X}^2)&=\displaystyle E((\frac{1}{n} \sum\limits_{i=1}^n \xi_i)^2)
=\displaystyle E((\frac{1}{n} \sum\limits_{i=1}^n (\xi_i-\mu+\mu))^2) \\
&=\displaystyle \frac{1}{n}\sigma^2+\mu^2, \\
\end{array}
\end{equation*}

\begin{equation*}
\begin{array}{ll}
\displaystyle E(\widetilde{X}^2)&=\displaystyle E((\frac{1}{n} \sum\limits_{i\in J\setminus J_{k^*}}\xi_i)^2)
=\displaystyle E((\frac{1}{n}\sum\limits_{i\in J\setminus J_{k^*}} (\xi_i-\mu+\mu))^2)\\
&=\displaystyle (\frac{1}{n}-\frac{k^*}{n^2})\sigma^2+(1-\frac{k^*}{n})^2 \mu^2, \\
\end{array}
\end{equation*}

then,
\begin{equation*}
\begin{array}{ll}
\displaystyle E\widetilde{S}_{1}^2 &
=\displaystyle E(\displaystyle \frac{1}{n}\sum\limits_{i=1}^n \xi_{i}^2-2\widetilde{X}\overline{X}+\widetilde{X}^2)\\
&=\displaystyle (1-\frac{1}{n}+\frac{k^{*}}{n^2})\sigma^2+\frac{{k^*}^2}{n^2}\mu^2.\\
\end{array}
\end{equation*}

\begin{equation*}
\begin{array}{ll}
\displaystyle E\widetilde{S}_{2}^2&=\displaystyle E(\frac{1}{n}\sum\limits_{i\in J\setminus J_{k^*}} \xi_{i}^2-2\widetilde{X}\overline{X}+(1-\frac{k^{*}}{n})\overline{X}^2),\\
&=\displaystyle (1-\frac{1}{n}-\frac{k^{*}}{n}+\frac{k^{*}}{n^2})\sigma^2.\\
\end{array}
\end{equation*}

\begin{equation*}
\begin{array}{ll}
\displaystyle E\widetilde{S}_{3}^2&=E(\displaystyle \frac{1}{n}\sum\limits_{i\in J\setminus J_{k^*}} \xi_i^2-(1+\frac{k^{*}}{n})\widetilde{X}^2)\\
 &=\displaystyle (1-\frac{1}{n}-\frac{k^{*}}{n}+\frac{k^{*}}{n^3})\sigma^2+(1-\frac{k^*}{n})\frac{k^*}{n^2}\mu^2.\\
\end{array}
\end{equation*}

Hence, it completes the proof.

\hfill $\Box$

\textbf{Corollary 1} Suppose $\xi_{1}, \xi_{2}, \cdots, \xi_{n}, \cdots $ are i.i.d random  variables of  $\xi$,  with $E\xi=\mu$, $D\xi=\sigma^2$. Then,
\begin{equation}
\begin{array}{ll}
\displaystyle E\widetilde{S}_{1}^2 > ES^2,\\
\displaystyle E\widetilde{S}_{2}^2 < ES^2, \\
\displaystyle E\widetilde{S}_{3}^2 \left\{\begin{array}{cl}
                   < ES^2,&\mbox{if}\quad \displaystyle \mu=0.\\
                   \leq ES^2,&\mbox{if}\quad \displaystyle \mu\neq 0, k^*\geq n-(n^2-1)\frac{\sigma^2}{\mu^2}.\\
                   > ES^2,&\mbox{if} \quad \displaystyle \mu\neq 0, k^*< n-(n^2-1)\frac{\sigma^2}{\mu^2}.\\
                                     \end{array} \right.
\end{array}
\end{equation}

\textbf{Proof}
Since,
\begin{equation*}
\begin{array}{ll}
\displaystyle E\widetilde{S}_{1}^2
&=\displaystyle (1-\frac{1}{n}+\frac{k^{*}}{n^2})\sigma^2+\frac{{k^*}^2}{n^2}\mu^2.\\
&=\displaystyle (1-\frac{1}{n})\sigma^2+(\frac{k^{*}}{n^2}\sigma^2+\frac{{k^*}^2}{n^2}\mu^2).\\
&> ES^2.
\end{array}
\end{equation*}

\begin{equation*}
\begin{array}{ll}
\displaystyle E\widetilde{S}_{2}^2 &=\displaystyle (1-\frac{1}{n}-\frac{k^{*}}{n}+\frac{k^{*}}{n^2})\sigma^2.\\
&=\displaystyle (1-\frac{1}{n})\sigma^2 -\frac{k^{*}}{n}(1-\frac{1}{n})\sigma^2.\\
&< ES^2.
\end{array}
\end{equation*}

\begin{equation*}
\begin{array}{ll}
\displaystyle E\widetilde{S}_{3}^2
&=\displaystyle (1-\frac{1}{n}-\frac{k^{*}}{n}+\frac{k^{*}}{n^3})\sigma^2+(1-\frac{k^*}{n})\frac{k^*}{n^2}\mu^2.\\
\end{array}
\end{equation*}

If $\mu=0$,
\begin{equation*}
\begin{array}{ll}
\displaystyle E\widetilde{S}_{3}^2
&=\displaystyle (1-\frac{1}{n}-\frac{k^{*}}{n}+\frac{k^{*}}{n^3})\sigma^2+(1-\frac{k^*}{n})\frac{k^*}{n^2}\mu^2.\\
&=\displaystyle (1-\frac{1}{n})\sigma^2-\frac{k^{*}}{n}(1-\frac{1}{n^2})\sigma^2.\\
&< ES^2.
\end{array}
\end{equation*}

If $\mu\neq 0$, then
\begin{equation*}
\begin{array}{ll}
\displaystyle E\widetilde{S}_{3}^2
&=\displaystyle (1-\frac{1}{n}-\frac{k^{*}}{n}+\frac{k^{*}}{n^3})\sigma^2+(1-\frac{k^*}{n})\frac{k^*}{n^2}\mu^2.\\
&=\displaystyle (1-\frac{1}{n})\sigma^2-[\frac{k^{*}}{n}(1-\frac{1}{n^2})\sigma^2-(1-\frac{k^*}{n})\frac{k^*}{n^2}\mu^2].\\
\end{array}
\end{equation*}

If $\displaystyle \frac{k^{*}}{n}(1-\frac{1}{n^2})\sigma^2-(1-\frac{k^*}{n})\frac{k^*}{n^2}\mu^2\geq 0$, then $\displaystyle E\widetilde{S}_{3}^2 \leq ES^2$.

That is
\begin{equation*}
\displaystyle \frac{k^{*}}{n}(1-\frac{1}{n^2})\sigma^2\geq (1-\frac{k^*}{n})\frac{k^*}{n^2}\mu^2,
\end{equation*}

which equals to
\begin{equation*}
\displaystyle k^*\geq n-(n^2-1)\frac{\sigma^2}{\mu^2}.
\end{equation*}

Therefore, If $\displaystyle k^*< n-(n^2-1)\frac{\sigma^2}{\mu^2}$, then $\displaystyle E\widetilde{S}_{3}^2 > ES^2$.

Hence, it completes the proof.

\hfill $\Box$

 From above discussion, $\displaystyle \frac{1}{n} \sum\limits_{i\in J\setminus J_{k^*}}\xi_i$ is asymptotic biased estimators of $E(\xi)=\mu$ according to deleting items WLLN. And $\widetilde{S}_{1}^2$, $\widetilde{S}_{2}^2$ and $\widetilde{S}_{3}^2$ are the asymptotic biased estimators of variance $D(\xi)=\sigma^2 $, with different expectations from $S^2$.

\section{Conclusion}
\label{}

We develop the deleting items limit theories of random variable sequence by substituting the partial sum with deleting items sum of random variables. The deleting items limit theorems extend the classical WLLN, SLLN and CLT, and provide asymptotic estimator for sample expectation and variance. Our research provide many probability and statistical conclusions for theoretical and real-world application, especially in large data analysis, statistical learning, pattern recognition, etc. Here, we only address main limit theorems in probability theory, the future work will focus on more limit theory conclusions in probability and statistics and extend them to deleting items general style.

\section*{References}
\label{}

\end{document}